\documentclass[preprint,authoryear,12pt]{elsarticle}
\usepackage{amsfonts,amssymb,amsmath,amsopn}
\usepackage{graphics}
\usepackage{rotating}
\usepackage{bigints}
\usepackage{booktabs}
\usepackage{subfigure}
\usepackage{url}
\usepackage{amssymb}
\usepackage{color}
\usepackage{multirow}
\usepackage{diagbox}
\usepackage{enumerate}
\usepackage{float}
\usepackage{epstopdf}
\usepackage[colorlinks=true]{hyperref}
\newcommand{\ignore}[1]{}
\newtheorem{asu}{{\sc Assumption}}

\newtheorem{thm}{Theorem}
\newtheorem{lem}{Lemma}
\newdefinition{rmk}{Remark}
\usepackage[margin=1in]{geometry}
\newproof{pf}{Proof}
\newproof{pot}{Proof of Theorem \ref{thm2}}

\usepackage{amssymb}
\usepackage{booktabs}

\allowdisplaybreaks
\def\dfrac{\displaystyle\frac}

\newcommand{\red}{}
\journal{\null}
\begin{document}
	\begin{frontmatter}
		\title{
			Edgeworth corrections for spot volatility estimator
		}
		\author[UM]{Lidan HE}
		\author[NUS]{Qiang LIU\corref{cor1}}
		\author[UM]{Zhi LIU}
		\address[UM]{Department of Mathematics, University of Macau}
		\address[NUS]{Department of Mathematics, National University of Singapore}
		\cortext[cor1]{Corresponding author. Email: matliuq@nus.edu.sg} %Tel.: +853-88224494. }
		\cortext[cor2]{Qiang LIU's research is supported by MOE-AcRF grant R-146-000-258-114 of Singapore. Zhi LIU's research is supported by the Science and Technology Development Fund, Macau SAR (File no. 202/2017/A3) and NSFC(No. 11971507).  }
		\begin{abstract}
			We develop Edgeworth expansion theory for spot volatility estimator under general assumptions on the log-price process that allow for drift and leverage effect. The result is based on further estimation of skewness and kurtosis, when compared with existing second order asymptotic normality result. Thus our theory can provide with a refinement result for the finite sample distribution of spot volatility. We also construct feasible confidence intervals (one-sided and two-sided) for spot volatility by using Edgeworth expansion. The Monte Carlo simulation study we conduct shows that the intervals based on Edgeworth expansion perform better than the conventional intervals based on normal approximation, which justifies \red{the} correctness of our theoretical conclusion.  \\~\\
		\end{abstract}
	\end{frontmatter}
	\section{Introduction}
	The fast development of computer technology and its wide application in financial market have made high frequency data to be increasingly available. And its research on both statistics and econometrics has been experiencing a great growth over the last several decades. Volatility of an asset quantifies the strength of its \red{fluctuation} over time. It plays a pivotal role in the fields of asset and derivations pricing, portfolio selection, risk management, and hedging, etc. 
	
	Recently, spot volatility estimation by using high frequency data has been received substantial attention, since it enables one to determine the variation of an asset at any given time. From a theoretical point of view, if we model the latent price of an asset as a continuous semi-martingale, spot volatility is just the coefficient of diffusion part, namely the conditional variance of the price. By rolling and blocking sampling filters,  \citet{FN1996} estimated spot volatility from high frequency data for the first time, and proved a pointwise asymptotic normality for rolling regression estimators. In \citet{FW2008}, the researchers proposed a kernel type estimator for spot volatility and established its explicit asymptotic distribution, when the price and volatility processes of an asset are modeled by bivariate diffusion \red{processes}. More literatures on kernel smoothing \red{for} the estimation of spot volatility, where microstructure noise or jumps may be accommodated, can be referred to \citet
{R2008}, \citet{D2010}, \citet{YP2014}, \citet{YFLZZ2014}, \citet{LLL2018} and references therein. 
	
	%\citet{FN1996} firstly showed that spot volatility can be estimated from high frequency data by rolling and blocking sampling filters and established a pointwise asymptotic normality for rolling regression estimators. \citet{FW2008} employed a bivariate diffusion to model the price and volatility of an asset and investigates kernel type estimators of spot volatility based on high-frequency return data and established both point-wise and global asymptotic distributions for the estimators. (\citet{ZH2014}) constructed a spot volatility estimator for high-frequency financial data which contain market micro-structure noise and proved consistency and derived the asymptotic distribution of the spot volatility estimator. (\citet{LLL2018}) proposed a kernel estimator for the spot volatility of a semi-martingale at a given time point by using high frequency data, where the underlying process accommodates a jump part of infinite variation. 
	 
	 Based on the asymptotic normality of the estimator of spot volatility, statistical inference on volatility can be made. More precisely, confidence intervals for spot volatility can be constructed. In this paper, our main motivation is to improve upon the existing asymptotic mixed normal approximation for the kernel estimator. Our theory is built upon general continuous semi-martingale assumption where a correlational relationship between the price and volatility processes, namely leverage effect in finance, is considered. 
	 
	 Edgeworth expansion is a power series result for the asymptotic distribution of an estimator that incorporates all moment information(see \citet{H1992} for a complete introduction). Thus, it can correct the asymptotic normal approximation by including the estimation of high order moments such as skewness and kurtosis. Recently, it has been applied to the estimation of volatility for correcting its performance in small samples. The Edgeworth expansion for realized volatility, which estimates the integrated volatility, was pioneeringly given in \citet{GM2009}. Their result was based on the assumption that the volatility process is independent of the price process, namely the leverage effect was ruled out, and the drift term \red{should not} be involved. By using \red{the aforementioned} conclusion, \citet{GM2008} discussed how confidence intervals could be constructed to correct normal approximation for realized volatility, and conducted some Monte \red{Carlo} simulation studies to validate their conclusion. \citet{ZMA2011} even \red{considered} the presence of microstructure noise when deriving Edgeworth expansions for realized volatility and other microstructure noise robust estimators. \citet{HV2016} established a full formal validity of Edgeworth expansions for realized volatility estimators given in above references. In this paper, we develop the theory of Edgeworth expansion for spot volatility estimator, and use it to construct corrected confidence intervals which refine conventional confidence intervals based on normal approximation. 
	 
	 %Their research has made enormous contributions to the financial market and has some limitations at the same time. They just considered the first order asymptotic normal theories. We are plan to give a more accurate asymptotic distributions for the spot volatility estimator based on Edgeworth expansions. Edgeworth expansion method can be dated from Chebyshev's paper in 1890 and is now a century old. It has recently received a revival of interest owing to its usefulness for exploring properities of contemporary statistical methods. (\citet{GM2008}) considered constructing confidence intervals for integrated volatility using Edgeworth expansions. Inspired by them, our main motivation is to improve the quality of the asymptotic normality of spot volatility using these Edgeworth expansions. In this article, we focus on spot volatility and ask whether we can improve upon the existing first-order asymptotic theory by relying on Edgeworth expansions in the absence of microstructure noise. Edgeworth expansions correct the first order asymptotic confidence interval by providing the explicit second order confidence interval. Using these approach is quiet important in model analysis and inference.   

    The paper is organized as follows. In Section \ref{sec:setup}, we give out the theoretical set up of our model and related assumptions. We simply review the spot volatility estimator of kernel type and develop its Edgeworth expansion in Section \ref{sec:main}, where the corrected confidence intervals are also constructed. In Section \ref{sec:simu}, some Monte Carlo simulation studies are conducted for evaluating the finite sample performance of our proposed corrected confidence intervals. Section \ref{sec:conc} concludes our paper. The theoretical proofs are deferred to Appendix part.

	\section{Setup}\label{sec:setup}
	Under the assumption of arbitrage-free and frictionless market, the logarithmic price of an asset $\{X_t \}_{t \in [0,T]}$ is necessarily to be modeled as a semi-martingale process (\citet{DS1994}). In this paper, we assume $\{X_t \}_{t \in [0,T]}$ is a continuous It$\hat{\text{o}}$ semi-martingale without the presence of jumps. It is a fundamental case that is most widely used in econometrics literatures. Under the continuous setting, the underlying data generating process $X_t$ defined on a filtered probability space $(\Omega, \mathcal{F}, \{\mathcal{F}_t\}_{t \in [0,T]}, \mathcal{P})$ is driven by 
	\begin{eqnarray}\label{model:con}
	dX_t=b_t dt+\sigma_t dB_t,\ t\in [0,T] ,
	\end{eqnarray}
	where $B$ is a standard Brownian motion, $b$ and $\sigma$ are adapted and locally bounded c$\acute{\text{a}}$dl$\acute{\text{a}}$g processes. To  guarantee the existence and uniqueness of the solution for the stochastic differential equation (\ref{model:con}), we assume the following Lipschitz continuity conditions are satisfied for the volatility process $\sigma$.
    \begin{asu}\label{assu:vol}
	For $s,t \in [0,T]$, there exists a constant $C$ and $0 < \alpha < 1$ such that 
		\begin{align*}
		\mathbf{E}[(\sigma_s - \sigma_t)^2] \leq C|s-t|^{\alpha}.
		\end{align*}
	Moreover, $\sigma_{t}^2$ is bounded away from 0, that is, there exists a constant $c$ such that $\sigma_{t}^2 > c > 0$.
	\end{asu}
We note that this is a rather general assumption and is widely used in many other literatures such as \citet{JT2014}, \citet{LLL2018}, \citet{YP2014}, etc. Possible models of $\sigma$ that satisfy the above assumption can be 
\begin{eqnarray}\label{model:vol}
d\sigma_t=b^{\sigma}_t dt+\sigma'_t \mathrm{d}B_t + \sigma''_t \mathrm{d}W_t,\ t\in [0,T] ,
\end{eqnarray}
where $W$ is another standard Brownian motion independent of $B$, and $b^{\sigma}, \sigma', \sigma''$ are adapted and locally bounded c$\acute{\text{a}}$dl$\acute{\text{a}}$g processes. In this case, assumption (\ref{assu:vol}) can be satisfied by taking $0 < \alpha \leq 1/2$. Further, the presence of jumps can also be involved in this model, which shall not violate the assumption. Interested readers can refer to \citet{JT2014} for the explicit form. We also note that the common driving process $B$ between the price process (\ref{model:con}) and the volatility process (\ref{model:vol}) depicts their correlated relationship, which is called leverage effect in finance. While in \citet{GM2009} and \citet{GM2008}, independent structure between $X$ and $\sigma$ is required for them to derive Edgeworth expansions for realized volatility. In this sense, our model is a general extension to their one, based on which the Edgeworth expansion for the spot volatility $\sigma_{\tau}^2$ at time $\tau$ is developed.

In practice, the whole realization path of $\{X_t \}$ for $t \in [0,T]$ is not achievable, and the price data are recorded at some finite time points. Without loss of generality, we assume the observations are obtained at fixed time points that are equally distributed within $[0,T]$, that is $\{0,\Delta_n,2\Delta_n,\cdots,n\Delta_n\}$ with $\Delta_n = \dfrac{T}{n}$. As $n$ tends to infinity, the length of time span for continuously observed data $\Delta_n$ shrinks, and it results in the so-called high frequency data. In what follows, our whole theory shall based on such an infill setting by taking $n\rightarrow \infty$. We define the shorthand $\Delta_i^nX:=X_{i\Delta_n}-X_{(i-1)\Delta_n}$ for $i=1,...,n$. 

\section{Main results}\label{sec:main}
\subsection{Spot volatility estimator}
In this paper, we are interested in estimating the spot volatility $\sigma^2_{\tau}$ at a given time $\tau \in [0,T]$. One of the most often used technique is by plugging in a kernel function into an estimator of the integrated volatility $\int_{0}^{T}\sigma^2_tdt$ and then letting the bandwidth parameter tends to 0 (see, e.g. \citet{FW2008}, \citet {R2008}, \citet{D2010}, \citet{YP2014}, \citet{YFLZZ2014}, \citet{LLL2018}).
\red{
Namely, the kernelized estimator of $\sigma_{\tau}^2$ when realized volatility in \citet{BNS2002a} is applied can be written as 
	\begin{align}\label{genker}
	 \widehat{\sigma^2_{\tau}}^{ker}= \Delta_n\sum_{i=1}^{n} K_h(i\Delta_n-\tau) (\Delta_i^n X)^2,
	\end{align}
where $h$ is the bandwidth parameter, $K_h(x) = K(x/h)/h $ with $K(x)$ being the kernel function defined on bounded interval $[a,b]$. We also assume that $K(x)$ is nonnegative and continuously differentiable with 
\begin{equation}
\int_a^bK^2(x)dx < \infty,~ \int_a^b K(x)dx =1.
\end{equation}
In this paper, we consider the specific kernel function of $K(x) = 1_{ \{0 \leq x < 1\} }$ for clarity:
	\begin{align}\label{uniker}
	\widehat{\sigma^2_{\tau}}= \dfrac{1}{k_n\Delta_n}\sum_{i=\lfloor \tau/\Delta_n  \rfloor + 1}^{\lfloor \tau /\Delta_n \rfloor +k_n } (\Delta_i^n X)^2,
	\end{align}
where $k_n := \lfloor h/\Delta_n \rfloor$ is the number of intraday returns that are close to time $\tau$ and approximately used for quantifying the variation of price process $X$ at that time. We note that the asymptotic properties of $\widehat{\sigma^2_{\tau}}$ can be generally extended for $\widehat{\sigma^2_{\tau}}^{ker}$ by lettting $1/k_n$ in (\ref{uniker}) to be $\frac{\Delta_n}{h} K( (i\Delta_n - \tau)/h )$ in (\ref{genker}). We see that for different kernel functions, different weights are used for the increments, which lead to possible different asymptotic variances and higher order moments for our use in this paper. This can be seen from (12) and (16) in \citet{LLL2018}, and uniform, Epanechnikov, quartic, triweight kernel functions are discussed there for illustration.
}

According to the asymptotic results given in the aforementioned existing literatures for the kernel version of the spot volatility estimator, we have  
\begin{align*}
\sqrt{k_n}(\widehat{\sigma^2_{\tau}} - \sigma_{\tau}^2) \rightarrow^{st} \mathcal{N}(0,2\sigma_{\tau}^4), \ \text{as} \ k_n\rightarrow \infty, \ k_n\Delta_n \rightarrow 0,
\end{align*}
where $\rightarrow^{st}$ means converging stably, which is a stronger result than convergence in distribution. Interested readers can refer to \citet{JS2003} for its \red{rigorous} definition and more detailed properties. And further, we have the following central limit theorem
	\begin{align}\label{con:norS}
	S(\tau,k_n):=\dfrac{\sqrt{k_n}(\widehat{\sigma^2_{\tau}} - \sigma_{\tau}^2)}{\sqrt{2}\sigma_{\tau}^2} \rightarrow^{d} \mathcal{N}(0,1),  \ \text{as} \ k_n\rightarrow \infty, \ k_n\Delta_n \rightarrow 0.
	\end{align}
   The above result is not feasible for \red{inferring} the information of $\sigma_{\tau}^2$ in practice since the denominator term of the statistic $S(\tau,k_n)$ relies on the underlying spot volatility $\sigma_{\tau}^2$. Since $\widehat{\sigma^2_{\tau}}$ can also be used to estimate $\sigma_{\tau}^2$ consistently, it gives us the following feasible version of second order asymptotic result:
	\begin{align}\label{con:norT}
	T(\tau,k_n):= \dfrac{\sqrt{k_n}(\widehat{\sigma^2_{\tau}} - \sigma_{\tau}^2)}{\sqrt{2}\widehat{\sigma^2_{\tau}}} \rightarrow^{d} \mathcal{N}(0,1),  \ \text{as} \ k_n\rightarrow \infty, \ k_n\Delta_n \rightarrow 0.
	\end{align}
	
	With the asymptotic distribution conclusions (\ref{con:norS}) and (\ref{con:norT}), statistical inference with respect to $\sigma_{\tau}^2$ turns to constructing confidence intervals for $\sigma_{\tau}^2$. In the proceeding, we will show how Edgeworth expansions can be derived for the statistics $S(\tau,k_n)$ and $T(\tau,k_n)$, based on which more accurate confidence interval results can be given. 

	\subsection{Edgeworth expansions for spot volatility estimator} 
	Let $k_j[S(\tau,k_n)], k_j[T(\tau,k_n)]$ denote the $j$-th order cumulant of $S(\tau,k_n)$ and $T(\tau,k_n)$. The Edgeworth expansions for $S(\tau,k_n)$ and $T(\tau,k_n)$ depend on their cumulants. The following lemma gives out the first fourth cumulants of the two statistics. 
	\begin{lem}\label{lem:cum}
	Under assumption \ref{assu:vol} and conditional on $\sigma_{\tau}$, we have 
	\begin{align*}
	k_1[S(\tau,k_n)]& =0 + O_p\big( k_n^{\alpha+1/2} \Delta_n^{\alpha} \big) , \quad k_2[S(\tau,k_n)] = 1+ O_p\big( k_n^{\alpha+1/2} \Delta_n^{\alpha} \big),  \\
	k_3[S(\tau,k_n)] &= \dfrac{2\sqrt{2}}{\sqrt{k_n}} + O_p\big( k_n^{\alpha+1/2} \Delta_n^{\alpha} \big), \quad k_4[S(\tau,k_n)]  = \dfrac{12}{k_n} +  O_p\big( k_n^{\alpha+1/2} \Delta_n^{\alpha} \big),
	\end{align*} 
   and further,
	\begin{align*}
	k_1[T(\tau,k_n)]& = \frac{-\sqrt{2}}{\sqrt{k_n}}+ O_p\big( k_n^{\alpha+\frac{1}{2}} \Delta_n^{\alpha} \big) +O_p\big(k_n^{-\frac{3}{2}}\big), \\
	k_2[T(\tau,k_n)]&= 1+\frac{8}{k_n} + O_p\big( k_n^{\alpha+\frac{1}{2}} \Delta_n^{\alpha} \big) +O_p\big(k_n^{-\frac{3}{2}}\big), \\
	k_3[T(\tau,k_n)] &= \dfrac{-4\sqrt{2}}{\sqrt{k_n}} + O_p\big( k_n^{\alpha+\frac{1}{2}} \Delta_n^{\alpha} \big) +O_p\big(k_n^{-\frac{3}{2}}\big),  k_4[T(\tau,k_n)] = \dfrac{60}{k_n} + O_p\big( k_n^{\alpha+\frac{1}{2}} \Delta_n^{\alpha} \big) +O_p\big(k_n^{-\frac{3}{2}}\big).
	\end{align*} 
\end{lem}

Now, we are ready to give the Edgeworth expansions of $S(\tau,k_n)$ and $T(\tau,k_n)$. 
\begin{thm}\label{thm:edge}
Under assumption \ref{assu:vol} and conditional on $\sigma_{\tau}$, if $k_n \rightarrow \infty$ and $k_n^{\alpha+3/2} \Delta_n^{\alpha} \rightarrow 0$, then we have the following second order Edgeworth expansions for $S(\tau,k_n)$ and  $T(\tau,k_n)$ \red{for} any given $x\in \mathbb{R}$:
\begin{align}\label{thm:Sedge}
\mathbf{P}(S(\tau,k_n) \leq x)&=\Phi(x)+\dfrac{1}{\sqrt{k_n}}p_1(x)\phi{(x)}+\dfrac{1}{k_n}p_2(x)\phi{(x)} + o(\dfrac{1}{k_n}),\\  \label{thm:Tedge}
\mathbf{P}(T(\tau,k_n) \leq x) &=\Phi(x)+\dfrac{1}{\sqrt{k_n}}q_1(x)\phi{(x)}+\dfrac{1}{k_n}q_2(x)\phi{(x)} + o(\dfrac{1}{k_n}),
\end{align}
with 
\begin{align*}
&p_1(x)=-\dfrac{\sqrt{2}}{3}(x^2-1), \ p_2(x)=-\dfrac{1}{2}H_{3}(x)-\dfrac{1}{9}H_{5} (x), \\
& q_1(x)=\sqrt{2}+\dfrac{2\sqrt{2}}{3}(x^2-1), \ q_2(x)=-5H_{1}(x)-\dfrac{23}{6}H_{3}(x)+\dfrac{4}{9}H_{5}(x),
\end{align*}
where $\Phi(\cdot)$ and $\phi(\cdot)$ are the standard normal cumulative and partial distribution functions respectively, $H_{i}$ denotes the $i$-th order Hermite polynomials with $H_{1}(x)=x,$\  $H_{3}(x)=x(x^2-3),$\  $H_{5}(x)=x(x^4-10x^2+15)$.
\end{thm}
\red{
	\begin{rmk}
	Considering the sample mean estimator of independent and identically distributed random variables, its tail probability is obtained by using the characteristic function. The Hermite polynomials are from the inverse Fourier-Stieltjes transform of the characteristic function of standard normal random variable and are orthogonal with respect to $\phi$. The detailed derivation can be found in Section 2.2 in \citet{H1992}. Thus, the conclusions of (\ref{thm:Sedge}) and (\ref{thm:Tedge}) above can be established if only the finite moment information of $S(\tau,k_n)$ and $T(\tau,k_n)$ can be approximated.
	\end{rmk}
}
\begin{rmk}
For the setting of the parameter $k_n$, one alternative is by taking $k_n = \lfloor  c\Delta_n^{-s} \rfloor$. In this case, the condition $k_n^{\alpha+3/2} \Delta_n^{\alpha} \rightarrow 0$ is equivalent to choosing $s$ with $ 0 < s < \dfrac{\alpha }{\alpha + 3/2}$. 
\end{rmk}
	\subsection{Corrected confidence intervals}
	
	In this section, we provide the confidence intervals for $\sigma_{\tau}^2$ based on the Edgeworth expansions in \red{the} last part. We will firstly describe the one-sided intervals, which are easier to understand. The discussion for \red{the} two-sided confidence interval follows. All of our discussions focus on intervals for  $\sigma_{\tau}^2$ based on the studentized statistic $T(\tau,k_n)$.    
	\subsubsection{One-sided confidence interval}
	Based on the asymptotic normality result (\ref{con:norT}),  we see that the conventional $95\%$ level one-sided confidence interval for $\sigma_{\tau}^2$ can be written as 
    \begin{align*}
	\mathcal{I}^{N-T,1}=(0, \widehat{\sigma^2_{\tau}}-\frac{\sqrt{2}\widehat{\sigma^2_{\tau}}z_{0.05}}{\sqrt{k_n}}),
	\end{align*} 
	where $z_{0.05}=-1.645$ is the $5\%$ quantile of standard normal distribution.
	By using the second order Edgeworth expansion result for $T(\tau,k_n)$ in (\ref{thm:Tedge}), the one-sided confidence interval has coverage probability equal to 
	\begin{align}\label{covp:NT1}
	\mathbf{P}(\sigma^2_{\tau} \in \mathcal{I}^{N-T,1})
	&=\mathbf{P}(T(\tau,k_n)\geq z_{0.05})=1-\mathbf{P}(T(\tau,k_n) < z_{0.05})\nonumber\\
	&=1-\Big[\Phi(z_{0.05})+\dfrac{\phi(z_{0.05})q_1(z_{0.05})}{\sqrt{k_n}}+o(\frac{1}{k_n})\Big]\nonumber\\
	&=0.95-\dfrac{\phi(z_{0.05})q_1(z_{0.05})}{\sqrt{k_n}}+o(\frac{1}{k_n}).
	\end{align}
	It's obvious that the error in coverage probability of $\mathcal{I}^{N-T,1}$ is of order $O(\dfrac{1}{\sqrt{k_n}})$. This inspires us to consider the following corrected one-sided confidence interval for $\sigma_{\tau}^2$:
	\begin{align*}
	\mathcal{I}^{E-T,1}=(0, \widehat{\sigma^2_{\tau}}-\frac{\sqrt{2}\widehat{\sigma^2_{\tau}}z_{0.05}}{\sqrt{k_n}}+\frac{\sqrt{2}\widehat{\sigma^2_{\tau}}q_1(z_{0.05})}{k_n}),
	\end{align*}
	where we recall that $q_1(x)$ is defined in Theorem \ref{thm:edge}. The above interval brings in a skewness correction term, that is $\frac{\sqrt{2}\widehat{\sigma^2_{\tau}}q_1(z_{0.05})}{k_n}$. Now, the coverage probability of $ I^{E-T,1} $ is 
	\begin{align}\label{covp:ET1}
	\begin{split}
	\mathbf{P}(\sigma^2_{\tau} \in \mathcal{I}^{E-T,1})
	&=\mathbf{P}\Big(T(\tau,k_n)\geq z_{0.05}-\frac{q_1(z_{0.05})}{\sqrt{k_n}}\Big)\\
    &=\Phi\Big(z_{0.05}-\frac{q_1(z_{0.05})}{\sqrt{k_n}}\Big)+\frac{q_1(z_{0.05}-\frac{q_1(z_{0.05})}{\sqrt{k_n}})}{\sqrt{k_n}}\phi\Big(z_{0.05}-\frac{q_1(z_{0.05})}{\sqrt{k_n}}\Big)+o\Big(\frac{1}{\sqrt{k_n}}\Big)\\
    & =0.95+O(\frac{1}{k_n}),
    \end{split}
	\end{align}
	which follows from arguments in Section 3.8 of \citet{H1992}. We see that the coverage probability error for $ I^{E-T,1}$ is of order $O(\frac{1}{k_n})$. Compared with the order of $O(\frac{1}{\sqrt{k_n}})$ for $\mathcal{I}^{N-T,1}$ based on the normal approximation, the corrected interval provides us with more exact result.
	
	\subsubsection{Two-sided confidence interval}
     Similarily as the one-sided corrected confidence interval for $\sigma^2_{\tau}$ by applying Edgeworth expansion, we can also develop \red{the} corresponding two-sided version. Following the discussion in the last part, by using the asymptotic normality result (\ref{con:norT}), the conventional $95\%$ level two-sided confidence interval for $\sigma_{\tau}^2$ is
	 \begin{align}
	 \mathcal{I}^{N-T,2}=(\widehat{\sigma^2_{\tau}}-\frac{\sqrt{2}\widehat{\sigma^2_{\tau}}z_{0.975}}{\sqrt{k_n}}, \widehat{\sigma^2_{\tau}}+\frac{\sqrt{2}\widehat{\sigma^2_{\tau}}z_{0.975}}{\sqrt{k_n}}),
	 \end{align}  
	 where $z_{0.975}=1.96$ is the $97.5\%$ quantile of standard normal distribution.
	 Its coverage probability is given by 
   \begin{align}\label{covp:NT2}
   \begin{split}
   P(\sigma^2_{\tau} \in \mathcal{I}^{N-T,2})
   &=P(|T(\tau,k_n)|\leq z_{0.975})\\
   &=2\Phi(z_{0.975})-1+2\frac{\phi(z_{0.975})q_2(z_{0.975})}{k_n}+o\Big(\frac{1}{k_n}\Big)\\
   &=0.95+2\frac{\phi(z_{0.975})q_2(z_{0.975})}{k_n}+o(\frac{1}{k_n}).
   \end{split}
   \end{align} 
    The above result is derived by using the second order Edgeworth expansion result for $T(\tau,k_n)$--(\ref{thm:Tedge}), together with the symmetry of $\Phi$, $\phi$, $q_1$ and $q_2$. 
      It can be seen that the error oder of coverage probability for $\mathcal{I}^{N-T,2}$ is $O(\frac{1}{k_n})$. 
      The corrected interval which contains a skewness and kurtosis correction term, and is based on the Edgeworth expansion of $T(\tau,k_n)$, is given as
	\begin{align*}
	\mathcal{I}^{E-T,2}=\Big(\widehat{\sigma^2_{\tau}}-\frac{\sqrt{2}\widehat{\sigma^2_{\tau}}z_{0.975}}{\sqrt{k_n}}+\frac{\sqrt{2}\widehat{\sigma^2_{\tau}}q_2(z_{0.975})}{k_n^{\frac{3}{2}}}, \widehat{\sigma^2_{\tau}}+\frac{\sqrt{2}\widehat{\sigma^2_{\tau}}z_{0.975}}{\sqrt{k_n}}-\frac{\sqrt{2}\widehat{\sigma^2_{\tau}}q_2(z_{0.975})}{k_n^{\frac{3}{2}}} \Big).
	\end{align*}
	 By similar proof as (\ref{covp:ET1}), we can show that the coverage probability of $ \mathcal{I}^{E-T,2}$ is
	\begin{align}\label{covp:ET2}
	\begin{split}
	P(\sigma^2_{\tau} \in \mathcal{I}^{E-T,2})
	&=P(|T(\tau,k_n)|\leq z_{0.05}-\frac{q_2(z_{0.975})}{k_n}) =0.95+O(\dfrac{1}{k_n^2}),
	\end{split}
	\end{align} 
	which implies that the coverage probability error order of  $I^{E-T,2}$ is $O(\dfrac{1}{k_n^2})$. Comparing the results (\ref{covp:NT2}) and (\ref{covp:ET2}) demonstrates us to what degree the two-sided confience interval is corrected by using the Edgewroth expansion derived. 
	
	Both the one-sided  and two-sided corrected confidence intervals have a smaller error order than the corresponding ones for normal approximation. Until now, we have provided the corrected confidence intervals for $\sigma^2_{\tau}$ based on the studentized statistic $T(\tau,k_n)$. In fact, similar results also hold for the normalized statistic $S(\tau,k_n)$. But since it is an infeasible statistic, we \red{do not} give a detailed discussion on it.
	
	\section{Simulation studies}\label{sec:simu}
	In this section, we conduct some Monte Carlo studies to evaluate the finite sample performance of the corrected intervals based on the Edgeworth expansion, namely $\mathcal{I}^{E-T,1}$ and $\mathcal{I}^{E-T,2}$. We also compare their performance with the one of respective asymptotic theory-based intervals $\mathcal{I}^{N-T,1}$ and $\mathcal{I}^{N-T,2}$. The simulation results show that the corrected intervals always outperform corresponding non-corrected versions under different settings, which verifies our theoretical analyses in \red{the} last section. 
	
	We consider two stochastic volatility models in our data generating process (\ref{model:con}). One of them is the following one factor stochastic volatility model
	\begin{align*}
	& \red{\text{Model}}~ \uppercase\expandafter{\romannumeral 1} : \quad \sigma_{t}= \exp(\beta_0+\beta_1 v_t) ,\ dv_{t}=\alpha v_t dt+dW_{t},
	\end{align*}   
	where $W$ is a standard Brownian motion independent of $B$; $\beta_0$, $\beta_1$ and $\alpha$ are constants. The other one is a two factor stochastic volatility model: 
	\begin{align*}
	\red{\text{Model}}~ \uppercase\expandafter{\romannumeral 2} :  \quad
	\sigma_{t}&= f(\beta_0+\beta_1 v_{1t}+\beta_2 v_{2t}) ,\\
	dv_{1t}&=\alpha_1 v_{1t} dt+dW_{1t}, \	dv_{2t}=\alpha_2 v_{2t} dt+(1+\phi v_{2t})dW_{2t},
	\end{align*} 
	where $W_{1}$, $W_{2}$ are mutually independent standard Brownian motions and they are also independent of $B$; $\beta_0, \beta_1, \beta_2, \alpha_1, \alpha_2, \phi$ are constants;  and the function $f(\cdot)$ is defined as 
	\begin{equation*}
	f(x)= \left\{
	\begin{aligned}
	&\exp(x), & \red{\text{if}}~ x\leq \log(1.5),\\
	&1.5\sqrt{1-\log(1.5)+x^2/\log(1.5)}, & \text{otherwise}.
	\end{aligned}
	\right.
	\end{equation*}
	For the parameters setting, we follow the ones in \citet{YP2014}, \citet{HT2005} and \citet{BNHLS2008a} with $\beta_1=0.125$, $\alpha=-0.025$,  $\beta_0=\beta_1/(2\alpha)$ \red{ for \text{Model}~ $\uppercase\expandafter{\romannumeral 1}$ }; $\beta_0=-1.2$, $\beta_1=0.04$, $\beta_2=1.5$, $\alpha_1=-0.0037$, $\alpha_2=-1.386$, $\phi= 0.25$ \red{ for \text{Model}~ $\uppercase\expandafter{\romannumeral 2}$ }. The initial value of above models both are $0.1$. And we consider the drift term in (\ref{model:con}) is $b_t \equiv 1$. For aforementioned models, a total number of 10000 paths are generated, and the estimation of $\sigma_{\tau}^2$ at $\tau = 0.3, 0.5, 0.7$ are considered. Different choices of $n$ as $780$, $1560$, $4680$, $7800$, $11700$, $23400$ are considered, and they correspond to``30-second", ``15-second", ``5-second", ``3-second",``2-second", ``1-second" interval returns. We set $k_n$ as $\lfloor cn^{1/4} \rfloor$ with $c$ equals $0.5$.
	
	\begin{table}[H]
		\small \caption{Coverage probabilities of normal $95\%$ confidence intervals for $\sigma^2_{\tau}$ in $\red{\text{Model}}~ \uppercase\expandafter{\romannumeral 1}$ }
		\label{Table1}
		\centering
		{\begin{tabular}{ccccccccccc}
				\hline
				 &  \quad & \multicolumn{2}{c}{$\tau = 0.3$} & & \multicolumn{2}{c}{$\tau = 0.5$} & & \multicolumn{2}{c}{$\tau = 0.7$} \\
				\hline
				$n$  & \quad   & $\mathcal{I}^{N-T,1}$& $\mathcal{I}^{E-T,1}$ & \quad & $\mathcal{I}^{N-T,1}$ & 
				$\mathcal{I}^{E-T,1}$ & \quad  & $\mathcal{I}^{N-T,1}$ & $\mathcal{I}^{E-T,1}$\\
				\hline
				780   &\quad & 79.96 & 87.38 &\quad& 79.99 & 87.72 &\quad & 79.01 & 87.38&\quad \\
				1560  &\quad & 85.10 & 91.52 &\quad& 85.68 & 92.26 &\quad & 85.41 & 91.52&\quad \\
				4680  &\quad & 88.01 & 93.45 &\quad& 88.59 & 93.95 &\quad & 88.77 & 93.45&\quad \\
				7800  &\quad & 88.50 & 93.56 &\quad& 87.54 & 92.98 &\quad & 88.64 & 93.56&\quad \\
				11700 &\quad & 91.29 & 95.32 &\quad& 90.64 & 94.84 &\quad & 91.18 & 95.32&\quad \\
				23400 &\quad & 93.03 & 96.52 &\quad& 92.64 & 96.10 &\quad & 92.44 & 96.52&\quad \\
				\hline
				\hline
				$n$  & \quad   & $\mathcal{I}^{N-T,2}$& $\mathcal{I}^{E-T,2}$ & \quad & $\mathcal{I}^{N-T,2}$ &$\mathcal{I}^{E-T,2}$ & \quad  & $\mathcal{I}^{N-T,2}$ & $\mathcal{I}^{E-T,2}$\\
				\hline
				780   &\quad & 80.56 & 82.94 &\quad& 81.30 & 83.82 &\quad & 81.47 & 84.02&\quad \\
				1560  &\quad & 87.16 & 89.17 &\quad& 87.17 & 89.57 &\quad & 86.83 & 88.92&\quad \\
				4680  &\quad & 90.15 & 91.84 &\quad& 89.68 & 91.61 &\quad & 90.41 & 92.49&\quad \\
				7800  &\quad & 89.40 & 91.24 &\quad& 90.20 & 92.09 &\quad & 89.86 & 91.61&\quad \\
				11700 &\quad & 92.36 & 94.04 &\quad& 92.43 & 93.73 &\quad & 92.51 & 94.15&\quad \\
				23400 &\quad & 94.62 & 95.84 &\quad& 93.57 & 95.02 &\quad & 94.16 & 95.41&\quad \\
				\hline
		\end{tabular}}
	\end{table}
	
	Tables \ref{Table1}-\ref{Table2} record the coverage probabilities of $\mathcal{I}^{N-T, 1}$, $\mathcal{I}^{E-T, 1}$, $\mathcal{I}^{N-T, 2}$ and $\mathcal{I}^{E-T, 2}$, when a standard normal coverage probability of $95\%$ is considered for the above two models. Similar \red{phenomena} are observed for these two different models. The degrees of undercoverage for the normal approximation based intervals are larger than the ones for corresponding Edgeworth corrected versions. We see that for relative lower frequency data, namely smaller value of $n$, the degree of undercoverage is larger. When the frequency is high enough, say $n=23400$, the coverage probabilities for the Edgeworth expansion corrected confidence intervals almost equal to $95\%$. In short, the correction eliminates the coverage distortions associated with the conventional confidence intervals with good effect.
	
	\begin{table}[H]
		\small \caption{Coverage probabilities of normal $95\%$ confidence intervals for $\sigma^2_{\tau}$ in $\red{\text{Model}}~ \uppercase\expandafter{\romannumeral 2}$ }
		\label{Table2}
        \centering
			{\begin{tabular}{ccccccccccc}
					\hline
					 &  \quad & \multicolumn{2}{c}{$\tau = 0.3$ } & & \multicolumn{2}{c}{$\tau = 0.5$} & & \multicolumn{2}{c}{$\tau = 0.7$} \\
					\hline
				$n$  & \quad   & $\mathcal{I}^{N-T,1}$& $\mathcal{I}^{E-T,1}$ & \quad & $\mathcal{I}^{N-T,1}$ &$\mathcal{I}^{E-T,1}$ & \quad  & $\mathcal{I}^{N-T,1}$ & $\mathcal{I}^{E-T,1}$\\
					\hline
					780   &\quad & 78.95 & 86.52 &\quad& 78.75 & 86.68 &\quad & 79.15 & 87.02&\quad \\
					1560  &\quad & 86.11 & 92.14 &\quad& 85.00 & 91.74 &\quad & 85.13 & 91.87&\quad \\
					4680  &\quad & 88.80 & 93.95 &\quad& 88.75 & 93.97 &\quad & 88.97 & 94.22&\quad \\
					7800  &\quad & 88.18 & 93.15 &\quad& 88.48 & 93.45 &\quad & 88.20 & 93.19&\quad \\
					11700 &\quad & 91.00 & 95.30 &\quad& 90.75 & 95.00 &\quad & 90.63 & 94.97&\quad \\
					23400 &\quad & 92.81 & 96.31 &\quad& 92.41 & 95.97 &\quad & 93.09 & 96.67&\quad \\
					\hline
					\hline
					$n$  & \quad   & $\mathcal{I}^{N-T,2}$& $\mathcal{I}^{E-T,2}$ & \quad & $\mathcal{I}^{N-T,2}$ &$\mathcal{I}^{E-T,2}$ & \quad  & $\mathcal{I}^{N-T,2}$ & $\mathcal{I}^{E-T,2}$\\
					\hline
				    780   &\quad & 80.79 & 83.36 &\quad& 81.02 & 83.47 &\quad & 80.77 & 83.34&\quad \\
				    1560  &\quad & 87.20 & 89.27 &\quad& 87.41 & 89.66 &\quad & 86.91 & 89.02&\quad \\
				    4680  &\quad & 90.45 & 92.22 &\quad& 90.22 & 92.27 &\quad & 90.19 & 92.12&\quad \\
				    7800  &\quad & 89.33 & 91.29 &\quad& 89.32 & 91.38 &\quad & 89.57 & 91.52&\quad \\
				    11700 &\quad & 92.36 & 93.98 &\quad& 92.81 & 94.46 &\quad & 92.31 & 93.92&\quad \\
				    23400 &\quad & 94.17 & 95.52 &\quad& 94.28 & 95.05 &\quad & 93.87 & 95.00& \quad \\				
					\hline
			\end{tabular}}
	\end{table}
	
	\section{Conclusion}\label{sec:conc}
	We derive Edgeworth expansion for \red{the} kernel type estimator of \red{the} spot volatility, which provides more exact result of asymptotic distribution than usual mixed normal distribution. Our theory is established in the presence of leverage effect, which \red{has not} been considered in other existing literatures on Edgeworth corrections for volatility estimators. By applying our theoretical conclusion, we give out corrections of the confidence intervals, one-sided or two-sided, with respect to the ones based on usual central limit theorem. In simulations, the superior finite sample performance of the corrected confidence intervals is observed, both for one-sided and two-sided versions.
	
	\section*{Appendix}  
	For simplicity of the proof procedure, we consider $b_s \equiv 0$ since the drift term $b$ has no effect on the estimation of volatility. And we define the following notations in advance:
	\begin{align*}
	&\widehat{\sigma^2_{\tau}}' = \dfrac{1}{k_n\Delta_n}\sum_{i=\lfloor \tau/\Delta_n  \rfloor + 1}^{\lfloor \tau /\Delta_n \rfloor +k_n } (\sigma_{\tau}\Delta_i^nB)^2, \quad R(\tau,k_n) = \dfrac{\sqrt{k_n}(\widehat{\sigma^2_{\tau}} - \widehat{\sigma^2_{\tau}}')}{\sqrt{2}\widehat{\sigma^2_{\tau}}}, \quad R'(\tau,k_n) = \dfrac{\sqrt{k_n}(\widehat{\sigma^2_{\tau}} - \widehat{\sigma^2_{\tau}}')}{\sqrt{2}\sigma^2_{\tau}}, \\ & M(\tau,k_n)=\dfrac{\sqrt{k_n}(\widehat{\sigma^2_{\tau}}'-\sigma_{\tau}^2)}{\sqrt{2}\sigma_{\tau}^2},\ U(\tau,k_n)=\dfrac{\sqrt{k_n}(\widehat{\sigma^2_{\tau}}-\sigma_{\tau}^2)}{\sigma_{\tau}^2}, \ Q(\tau,k_n) =  M(\tau,k_n)(1+\dfrac{1}{\sqrt{k_n}}U(\tau,k_n))^{-1}
	\end{align*}
	and observe that
	\begin{align*}
	S(\tau,k_n) &= \dfrac{\sqrt{k_n}(\widehat{\sigma^2_{\tau}} - \sigma_{\tau}^2)}{\sqrt{2}\sigma^2_{\tau}} = \dfrac{\sqrt{k_n}(\widehat{\sigma^2_{\tau}} - \widehat{\sigma^2_{\tau}}')}{\sqrt{2}\sigma^2_{\tau}} + \dfrac{\sqrt{k_n}(\widehat{\sigma^2_{\tau}}' - \sigma_{\tau}^2)}{\sqrt{2}\sigma_{\tau}^2} = R'(\tau,k_n)  + M(\tau,k_n), \\
	T(\tau,k_n) &= \dfrac{\sqrt{k_n}(\widehat{\sigma^2_{\tau}} - \sigma_{\tau}^2)}{\sqrt{2}\widehat{\sigma^2_{\tau}}} = \dfrac{\sqrt{k_n}(\widehat{\sigma^2_{\tau}} - \widehat{\sigma^2_{\tau}}')}{\sqrt{2}\widehat{\sigma^2_{\tau}}} + \dfrac{\sqrt{k_n}(\widehat{\sigma^2_{\tau}}' - \sigma_{\tau}^2)}{\sqrt{2}\sigma_{\tau}^2}\dfrac{\sigma_{\tau}^2}{\widehat{\sigma^2_{\tau}}}   = R(\tau,k_n) + Q(\tau,k_n).
	\end{align*}
	
\textbf{Proof of Lemma \ref{lem:cum}:} For $R'(\tau,k_n)$ and $R(\tau,k_n)$, under assumption \ref{model:con}, we have 
\begin{align*}
&\mathbf{E}[|R(\tau,k_n)|] = \mathbf{E}[| \dfrac{\sqrt{k_n}(\widehat{\sigma^2_{\tau}} - \widehat{\sigma^2_{\tau}}')}{\sqrt{2}\widehat{\sigma^2_{\tau}}} | ]\\
& \leq \sqrt{k_n} \dfrac{C}{k_n\Delta_n}\sum_{i=\lfloor \tau/\Delta_n  \rfloor + 1}^{\lfloor \tau /\Delta_n \rfloor +k_n } \mathbf{E}[| (\Delta_i^nX)^2 - (\sigma_{\tau}\Delta_i^nB)^2|] \\
& \leq \sqrt{k_n} \dfrac{C}{k_n\Delta_n} \sum_{i=\lfloor \tau/\Delta_n  \rfloor + 1}^{\lfloor \tau /\Delta_n \rfloor +k_n } (\mathbf{E}[ |\int_{(i-1)\Delta_n}^{i\Delta_n} (\sigma_s - \sigma_{\tau}) dB_s|^2])^{1/2} (\mathbf{E}[ |\int_{(i-1)\Delta_n}^{i\Delta_n} (\sigma_s + \sigma_{\tau}) dB_s|^2])^{1/2}\\
& \leq Ck_n^{\alpha + 1/2} \Delta_n^{\alpha}.
\end{align*}
The same result also holds for $R'(\tau,k_n)$ and can be similarily derived. Thus, we have $R(\tau,k_n) = O_p(k_n^{\alpha + 1/2} \Delta_n^{\alpha})$ and $R'(\tau,k_n) = O_p(k_n^{\alpha + 1/2} \Delta_n^{\alpha})$. 

For $Q(\tau,k_n) = M(\tau,k_n)(1+\dfrac{1}{\sqrt{k_n}}U(\tau,k_n))^{-1}$, according to the second order Taylor expansion of $f(x)=(1+x)^{-k}$ around $0$ for any fixed positive integer $k$, namely $f(x)=1-kx+\frac{k(k+1)}{2}x^2+O(x^3)$, together with the fact $M(\tau,k_n) = O_p(1)$ and $U(\tau,k_n)=O_p(1)$, we have 
\begin{align}\label{equ:Q}
Q(\tau,k_n)^k = M(\tau,k_n)^k
(1-k\frac{U(\tau,k_n)}{\sqrt{k_n}}+\frac{k(k+1)}{2}\frac{U(\tau,k_n)^2}{k_n}) + O_p(k_n^{-3/2}).
\end{align}
We note that condition on the information at time point $\tau$, $\sigma_{\tau}^2$ can be seen as a constant, and the following results hold
\begin{align*}
& \mathbf{E}[M(\tau,k_n)] =0, \quad \mathbf{E}[M(\tau,k_n)^2] =1, \quad  \mathbf{E}[M(\tau,k_n)^3] =\dfrac{2\sqrt{2}}{\sqrt{k_n}}, \quad \mathbf{E}[M(\tau,k_n)^4] =3+\dfrac{12}{k_n}, \\
& \mathbf{E}[M(\tau,k_n)^5] =\dfrac{20\sqrt{2}}{\sqrt{k_n}} + \dfrac{48\sqrt{2}}{k_n^{3/2}}, \quad \mathbf{E}[M(\tau,k_n)^6] =15+\dfrac{260}{k_n} + \dfrac{480}{k_n^2}, \\
&\mathbf{E}[M(\tau,k_n)U(\tau,k_n)] =\sqrt{2}+O_p\big( k_n^{\alpha + 1/2} \Delta_n^{\alpha} \big), \quad \mathbf{E}[M(\tau,k_n)^2U(\tau,k_n)]= \frac{4}{\sqrt{k_n}}+O_p\big( k_n^{\alpha + 1/2} \Delta_n^{\alpha} \big),\\
&  \mathbf{E}[M(\tau,k_n)^3U(\tau,k_n)] = 3\sqrt{2} + \dfrac{12\sqrt{2}}{k_n} + O_p\big( k_n^{\alpha + 1/2} \Delta_n^{\alpha} \big), \\
& \mathbf{E}[M(\tau,k_n)^4U(\tau,k_n)] =\dfrac{40}{\sqrt{k_n}} + \dfrac{96}{k_n^{3/2}} + O_p\big( k_n^{\alpha + 1/2} \Delta_n^{\alpha} \big), \\
& \mathbf{E}[M(\tau,k_n)U(\tau,k_n)^2]= \dfrac{4\sqrt{2}}{\sqrt{k_n}} + O_p\big( k_n^{\alpha + 1/2} \Delta_n^{\alpha} \big), \mathbf{E}[M(\tau,k_n)^2U(\tau,k_n)^2] = 6 + \dfrac{24}{k_n} + 
O_p\big( k_n^{\alpha + 1/2} \Delta_n^{\alpha} \big), \\
& \mathbf{E}[M(\tau,k_n)^3U(\tau,k_n)^2] = \dfrac{40\sqrt{2}}{\sqrt{k_n}} + \dfrac{96\sqrt{2}}{k_n^{3/2}} + O_p\big( k_n^{\alpha + 1/2} \Delta_n^{\alpha} \big),\\
& \mathbf{E}[M(\tau,k_n)^4U(\tau,k_n)^2] =30+\dfrac{520}{k_n} + \dfrac{960}{k_n^2} + O_p\big( k_n^{\alpha + 1/2} \Delta_n^{\alpha} \big).
\end{align*}
Furthermore, from (\ref{equ:Q}) we obtain 
\begin{align*}
&\mathbf{E}[Q(\tau,k_n)]=\frac{-\sqrt{2}}{\sqrt{k_n}}+ O_p\big( k_n^{\alpha} \Delta_n^{\alpha} \big) +O_p\big(k_n^{-\frac{3}{2}}\big), \mathbf{E}[Q(\tau,k_n)^2]=1+\frac{10}{k_n}+O_p\big( k_n^{\alpha} \Delta_n^{\alpha} \big) +O_p\big(k_n^{-\frac{3}{2}}\big),\\
&\mathbf{E}[Q(\tau,k_n)^3]=\frac{-7\sqrt{2}}{\sqrt{k_n}}+O_p\big( k_n^{\alpha} \Delta_n^{\alpha} \big) +O_p\big(k_n^{-\frac{3}{2}}\big), \mathbf{E}[Q(\tau,k_n)^4]=3+\frac{152}{k_n}+O_p\big( k_n^{\alpha} \Delta_n^{\alpha} \big) +O_p\big(k_n^{-\frac{3}{2}}\big).
\end{align*}

Since $T(\tau,k_n) = R(\tau,k_n) + Q(\tau,k_n)$ and the first four cumulants of $T(\tau,k_n)$ are given by (see, e.g., \citet{H1992}):
\begin{align*}
k_1(T(\tau,k_n))&=\mathbf{E}[T(\tau,k_n)], \quad k_2(T(\tau,k_n))=\mathbf{E}[T(\tau,k_n)^2]-[\mathbf{E}[T(\tau,k_n)]]^2, \\
 k_3(T(\tau,k_n))&=\mathbf{E}[T(\tau,k_n)^3]-3\mathbf{E}[T(\tau,k_n)^2]\mathbf{E}[T(\tau,k_n)]+2[\mathbf{E}[T(\tau,k_n)]]^3, \\
 k_4(T(\tau,k_n))&=\mathbf{E}[T(\tau,k_n)^4]-4\mathbf{E}[T(\tau,k_n)^3]\mathbf{E}[T(\tau,k_n)]-3[\mathbf{E}[T(\tau,k_n)]^2]^2  \\
& \quad +12\mathbf{E}[T(\tau,k_n)^2][\mathbf{E}[T(\tau,k_n)]]^2-6[\mathbf{E}[T(\tau,k_n)]]^4,
\end{align*}
we further have
\begin{align*}
k_1[T(\tau,k_n)]& = \frac{-\sqrt{2}}{\sqrt{k_n}}+ O_p\big( k_n^{\alpha+\frac{1}{2}} \Delta_n^{\alpha} \big) +O_p\big(k_n^{-\frac{3}{2}}\big), k_2[T(\tau,k_n)]= 1+\frac{8}{k_n} + O_p\big( k_n^{\alpha+\frac{1}{2}} \Delta_n^{\alpha} \big) +O_p\big(k_n^{-\frac{3}{2}}\big), \\
k_3[T(\tau,k_n)] &= \dfrac{-4\sqrt{2}}{\sqrt{k_n}} + O_p\big( k_n^{\alpha+\frac{1}{2}} \Delta_n^{\alpha} \big) +O_p\big(k_n^{-\frac{3}{2}}\big),  k_4[T(\tau,k_n)] = \dfrac{60}{k_n} +  O_p\big( k_n^{\alpha+\frac{1}{2}} \Delta_n^{\alpha} \big) +O_p\big(k_n^{-\frac{3}{2}}\big).
\end{align*} 
And since $S(\tau,k_n) = R'(\tau,k_n) + M(\tau,k_n)$, similarily we obtain
\begin{align*}
k_1[S(\tau,k_n)]& =0 + O_p\big( k_n^{\alpha+1/2} \Delta_n^{\alpha} \big), \quad  k_2[S(\tau,k_n)]= 1+ O_p\big( k_n^{\alpha+1/2} \Delta_n^{\alpha} \big),  \\
k_3[S(\tau,k_n)] &= \dfrac{2\sqrt{2}}{\sqrt{k_n}} + O_p\big( k_n^{\alpha+1/2} \Delta_n^{\alpha} \big), \quad k_4[S(\tau,k_n)]  = \dfrac{12}{k_n} +  O_p\big( k_n^{\alpha+1/2} \Delta_n^{\alpha} \big).
\end{align*}\hfill $\square$

\textbf{Proof of Theorem \ref{thm:edge}:} We observe that
\begin{align}\label{edg:M}
\mathbf{P}(M(\tau,k_n) & \leq x)=\Phi(x)+\dfrac{1}{\sqrt{k_n}}p_1(x)\phi{(x)}+\dfrac{1}{k_n}p_2(x)\phi{(x)} + o(\dfrac{1}{k_n}),\\ \label{edg:Q}
\mathbf{P}(Q(\tau,k_n) &\leq x)=\Phi(x)+\dfrac{1}{\sqrt{k_n}}q_1(x)\phi{(x)}+\dfrac{1}{k_n}q_2(x)\phi{(x)} + o(\dfrac{1}{k_n}),
\end{align}
which follow from (2.17) and Section 2.3 in \citet{H1992}, the condition $k_n^{\alpha+3/2} \Delta_n^{\alpha} \rightarrow 0$, and the following cumulants 
\begin{align*}
k_1[M(\tau,k_n)]& =0 + O_p\big( k_n^{\alpha} \Delta_n^{\alpha} \big), \quad  k_2[M(\tau,k_n)]= 1+ O_p\big( k_n^{\alpha} \Delta_n^{\alpha} \big),  \\
k_3[M(\tau,k_n)] &= \dfrac{2\sqrt{2}}{\sqrt{k_n}} + O_p\big( k_n^{\alpha} \Delta_n^{\alpha} \big), \quad k_4[M(\tau,k_n)]  = \dfrac{12}{k_n} +  O_p\big( k_n^{\alpha} \Delta_n^{\alpha} \big),
\end{align*} 
and 
\begin{align*}
k_1[Q(\tau,k_n)]& = \frac{-\sqrt{2}}{\sqrt{k_n}}+ O_p\big( k_n^{\alpha} \Delta_n^{\alpha} \big) +O_p\big(k_n^{-\frac{3}{2}}\big), \ k_2[Q(\tau,k_n)]= 1+\frac{8}{k_n} + O_p\big( k_n^{\alpha} \Delta_n^{\alpha} \big) +O_p\big(k_n^{-\frac{3}{2}}\big), \\
k_3[Q(\tau,k_n)] &= \dfrac{-4\sqrt{2}}{\sqrt{k_n}} + O_p\big( k_n^{\alpha} \Delta_n^{\alpha} \big) +O_p\big(k_n^{-\frac{3}{2}}\big),
\ k_4[Q(\tau,k_n)] = \dfrac{60}{k_n} +  O_p\big( k_n^{\alpha} \Delta_n^{\alpha} \big) +O_p\big(k_n^{-\frac{3}{2}}\big). 
\end{align*} 
The above cumulant results can be easily seen from the proof of Lemma \ref{lem:cum}. 

For any given $x \in \mathbb{R}$ and $h \geq 0$, we note that there exists a constant $C$ such that 
\begin{align*}
\mathbf{P}(M(\tau,k_n) \leq x+h) - \mathbf{P}(M(\tau,k_n) \leq x)  \leq Ch,
\end{align*}
since $\Phi(x), p_1(x)\phi{(x)}, p_2(x)\phi{(x)}$ in (\ref{edg:M}) and $q_1(x)\phi{(x)}, q_2(x)\phi{(x)}$ in (\ref{edg:Q}) are differentiable with continuous derivative. As shown in the proof of Lemma \ref{lem:cum}, we have $R(\tau,k_n) = O_p(k_n^{\alpha + 1/2} \Delta_n^{\alpha})$ and $R'(\tau,k_n) = O_p(k_n^{\alpha + 1/2} \Delta_n^{\alpha})$. Thus,
\begin{align*}
\mathbf{P}(S(\tau,k_n) \leq x) &= \mathbf{P}(R'(\tau,k_n)  + M(\tau,k_n) \leq x) \\
&= \mathbf{P}( M(\tau,k_n) \leq x + O_p(k_n^{\alpha + 1/2} \Delta_n^{\alpha}))= \mathbf{P}( M(\tau,k_n) \leq x ) + O_p(k_n^{\alpha + 1/2} \Delta_n^{\alpha}), \\
\mathbf{P}(T(\tau,k_n) \leq x) &= \mathbf{P}(R(\tau,k_n)  + Q(\tau,k_n) \leq x) \\
&= \mathbf{P}( Q(\tau,k_n) \leq x + O_p(k_n^{\alpha + 1/2} \Delta_n^{\alpha}))= \mathbf{P}( Q(\tau,k_n) \leq x ) + O_p(k_n^{\alpha + 1/2} \Delta_n^{\alpha}).
\end{align*}
Together with the condition $k_n^{\alpha+3/2} \Delta_n^{\alpha} \rightarrow 0$, we obtain the conclusions (\ref{thm:Sedge}) and (\ref{thm:Tedge}). 
\hfill $\square$

	\hspace{-0.25in}\section*{\textbf{References}}
	\bibliographystyle{model2-names}
	\bibliography{liu}
\end{document}